\documentclass[10pt,twoside]{article}       
\usepackage{graphicx}

\def\al{\alpha}

\def\z{\mathbf 0}

\def\ie{i.e.}
\def\eg{e.g.}

\def\etc{ etc.}

\def\al{\alpha}
\def\bc{\begin{center}}
\def\ec{\end{center}}
\def\eqn{\end{equation}\noindent}
\def\eqnr{\end{eqnarray}\noindent}
\def\beqr{\begin{eqnarray}}

\begin{document}
\setcounter{page} {1}


\begin{center}
\Large {\bf A Note On Topological Conjugacy For  Perpetual Points }\\
 {\large Awadhesh Prasad\footnote{\emph{E-mail Address}:awadhesh@physics.du.ac.in}}\\
\par \centering 
\normalsize  Department of Physics and Astrophysics, University of Delhi, Delhi 110007, India.
\\
\end{center}

\begin{quote}
\textbf{Abstract:} Recently a new class of critical points, termed as {\sl perpetual points}, where acceleration 
becomes zero but the velocity remains non-zero, is observed in nonlinear dynamical systems. In this work we
show  whether a transformation also maps the perpetual points to another system  or not.
 We establish  mathematically that a linearly transformed system is
topologicaly conjugate, and hence does map the perpetual points. However,
for a nonlinear transformation, various other possibilities are also discussed.
It is noticed that under a linear diffeomorphic transformation, perpetual points are mapped, and 
accordingly, eigenvalues  are preserved.  \\

\end{quote}

\section{Introduction}
\label{sec-1}

\noindent 

Consider a  general dynamical system specified by the equations  
\begin{eqnarray}
\dot{\mathbf{X}}=\mathbf{f}(\mathbf{X})
\label{eq:f}
\end{eqnarray}
\noindent
where $\mathbf{X}=(x_1,x_2,...,x_n)^T \in\Sigma$ is $n$-dimensional vector of dynamical variables and 
 $\mathbf{f}=(f_1(\mathbf{X}), f_2(\mathbf{X}), ...f_n(\mathbf{X}))^T$ specifies the
 evolution equations (velocity vector) of the system.
Here $T$ stands for transpose of the vector.  Acceleration of  the system can be   obtained 
by taking time   derivative of Eq. (\ref{eq:f}), viz.,
\begin{eqnarray}
\nonumber
\ddot{\mathbf{X}} &=&D_{\mathbf{X}^T}\mathbf{f}(\mathbf{X})\cdot\mathbf{f}(\mathbf{X})\\
\label{gen}
                 &=&\mathbf{F}(\mathbf{X})
\label{eq:FF}
\end{eqnarray}
\noindent
where $\mathbf{F}=D_{\mathbf{X}^T}\mathbf{f}(\mathbf{X})\cdot\mathbf{f}(\mathbf{X})$\footnote{Here 
$D_{\mathbf{X}}$ stands for derivative with respect to  ${\mathbf{X}}$. }
 is termed as acceleration vector.

The fixed points ($\mathbf{X}_{FP}$) of system, Eq. (\ref{eq:f}), are the ones where velocity
 and acceleration are simultaneously zero.
Using linear stability analysis the dynamics near the fixed points  is determined by the 
eigenvalues $\lambda=\{\lambda_1,\lambda_2,....,\lambda_n\}$
 of the Jacobian $D_{\mathbf{X}^T}\mathbf{f}(\mathbf{X})$.
Similarly, the perpetual points  ($\mathbf{X}_{PP}$) are the ones where acceleration is zero while velocity
 remains  finite (nonzero).
The eigenvalues, $\mu =\{\mu_1,\mu_2,....,\mu_n\}$ of $D_{\mathbf{X}^T}\mathbf{F}(\mathbf{X})$ determine whether
velocity is either extremum or inflection type.
The understanding of motion around perpetual points is studied recently \cite{ap} where connection between
 $\lambda$ and $\mu$ are
discussed along with its various applications \cite{ap,ap2}.

Note that an understanding  of the topological conjugation between different systems is important in almost
 all branch of sciences. This, in fact, will suggest as to  whether the  behavior of one system  is similar to that of 
the  other or not. If there exists conjugacy between the systems 
then the systems behavior of one can immediately be derived from the knowledge of the system behavior of conjugate system.
Usually, under certain conditions, such conjugate systems
can be transformed into each other using an appropriate transformation.
Suppose there is transformation, $\mathbf{h}(\mathbf{X})$, from variables $\mathbf{X}$ to $\mathbf{Y}$, \ie,
\begin{eqnarray}
{\mathbf{Y}}=\mathbf{h}(\mathbf{X})
\label{eq:h}
\end{eqnarray}
\noindent
where $\mathbf{Y}=(y_1,y_2,...,y_n)^T \in\Omega$ and $\mathbf{h}=(h_1(\mathbf{X}), h_2(\mathbf{X}),....h_n(\mathbf{X}))^T$, 
then the transformed  version of Eqs. (\ref{eq:f}) and  (\ref{eq:FF}) take the forms
\begin{eqnarray}
\dot{\mathbf{Y}} &= &\mathbf{g}(\mathbf{Y}) =\mathbf{g}(\mathbf{h}(\mathbf{X}))\\
\label{eq:g}
 \mbox{and~~~}\ddot{\mathbf{Y}} &= &\mathbf{G}(\mathbf{Y})=\mathbf{G}(\mathbf{h}(\mathbf{X})), 
\label{eq:GG}
\end{eqnarray}
\noindent
respectively. Here $\mathbf{g}(\mathbf{Y})$ and $\mathbf{G}(\mathbf{Y})$ are the velocity and the
 acceleration vectors of the transformed system.
 The dynamics near the fixed points is very important to understand the system behavior. 
Rigorous studies have been done to understand the  individual  as well as conjugate systems 
around the fixed points.  We recapitulate such works in Definition 1 and Theorem 1 mainly for the 
sake of  completeness.  The details can be found  in any standard books on dynamical systems \cite{arr,guck,kuz,meiss,wig}.
In the present work we show the topological conjugacy of  perpetual points between the systems,  
Eqs. (\ref{eq:FF}) and (\ref{eq:GG}).  The results are discussed and representative examples are presented
 in Sec. II.  
We summarize the results and highlight their importance in Sec. III

\section{Results and Discussions}\label{sec-2}

\noindent
{\bf Definition 1:}\footnote{For completeness the Definition 1 and Theorem 1 are reviewed  here. For details 
see: Refs. \cite{arr,guck,kuz,meiss,wig}. } The dynamical systems
$ \dot{\mathbf{X}} = \mathbf{f}(\mathbf{X})$ and $\dot{\mathbf{Y}} = \mathbf{g}(\mathbf{Y})$,
are topologically conjugate
if there exists a homeomorphism\footnote{$\mathbf{h}(\mathbf{X})$ should be one-to-one, onto, continuous,
 and has its continuous inverse also.}
 $\mathbf{h}: \Sigma\to\Omega$
for their corresponding flows $\phi_t{(\mathbf{X})}\in R\times\Sigma$ and $\psi_t{(\mathbf{Y})}\in R\times \Omega$, 
such that for each $\mathbf{X}\in \Sigma$ and $t\in R$

\begin{eqnarray}
\psi_t(\mathbf{\mathbf{h}(\mathbf{X})}) &=& \mathbf{h}(\phi_t(\mathbf{X})).
\label{eq:con}
\end{eqnarray}
Alternatively, one  can recast it as
\begin{eqnarray}
\psi \circ \mathbf{h} &=& \mathbf{h} \circ \phi\\
\mbox{or~~~~~~~~~~~} \psi  &=&  \mathbf{h}  \circ \phi  \circ\mathbf{h}^{-1}
\end{eqnarray}

\noindent
where the composition symbol, $\circ$, means $\psi \circ \mathbf{h}=\psi_t(\mathbf{\mathbf{h}(\mathbf{X})})$, \etc
~$\fbox{~}$\\

\noindent
{\bf Theorem 1:}$^3$ If the velocity vectors $\mathbf{f}(\mathbf{X}) $ and $\mathbf{g}(\mathbf{Y})$
 are topologically conjugate under
$\mathbf{h}(\mathbf{X})$ then fixed points of $\mathbf{f}(\mathbf{X})$ are mapped to the fixed points 
of $\mathbf{g}(\mathbf{Y})$.\\

\noindent
{\sl Proof:} From the Definition 1, the time derivative of Eq. (\ref{eq:con}) gives
\begin{eqnarray}
\dot{\psi}_t(\mathbf{h}(\mathbf{X})) &=& D_{\mathbf{X}^T}\mathbf{h}(\mathbf{X})\cdot \dot{\phi}_t(\mathbf{X})\\
\label{eq:fp0}
\mbox{\ie~~} \mathbf{g}(\mathbf{h}(\mathbf{X})) &=& D_{\mathbf{X}^T}\mathbf{h}(\mathbf{X})\cdot \mathbf{f}(\mathbf{X}).
\label{eq:fp}
\end{eqnarray} 
\noindent
Here, $ D_{\mathbf{X}^T}\mathbf{h}(\mathbf{X})$ is a Jacobian of size $n\times n$.
 Since the velocity vector $\mathbf{f} (\mathbf{X})=\z$  at the fixed points and hence Eq. (\ref{eq:fp}) follows
 $\mathbf{g}(\mathbf{Y})=\z$.
Note that $\mathbf{g}(\mathbf{Y})=\z$ is also possible for   $D_{\mathbf{X}^T}\mathbf{h}(\mathbf{X})=\z$ whether $\mathbf{f}(\mathbf{X})
 \ne \z$ or $\mathbf{f} (\mathbf{X})= \z$. 
This suggests  that for the case of nonlinear $\mathbf{h}(\mathbf{X})$, particularly when the possibility of
$D_{\mathbf{X}}\mathbf{h}(\mathbf{X})=\z$ exists,
 some new fixed points may be created\footnote{Due to nonlinear function, $\mathbf{h}(\mathbf{X})$, the order of polynomial of
$\mathbf{Y}$ in  $\mathbf{g} (\mathbf{Y})$ may be higher than that of the $\mathbf{X}$ in $\mathbf{f} (\mathbf{X})$.} 
in addition to that of the $\mathbf{f}(\mathbf{X}) =\z$. $\fbox{~}$
~\\

\noindent
 {\bf Theorem 2:} If the velocity vectors $\mathbf{f}(\mathbf{X})$ and $\mathbf{g}(\mathbf{Y})$ are topologically conjugate under
homeomorphic$^4$ {\sl linear} $\mathbf{h}(\mathbf{X})$ then the perpetual points of $\mathbf{f}(\mathbf{X})$ are 
mapped to the perpetual points of $\mathbf{g}(\mathbf{Y})$.\\

\noindent
{\sl Proof:} 
The  time derivative of Eq. (\ref{eq:fp0}) gives
\begin{eqnarray}
\ddot{\psi}_t(\mathbf{h}(\mathbf{X})) & = & D^2_{t\mathbf{X}^T}\mathbf{h}(\mathbf{X})\cdot\dot{\phi}_t(\mathbf{X})+
D_{\mathbf{X}^T}\mathbf{h}(\mathbf{X})\cdot\ddot{\phi}_t(\mathbf{X})\\
\mbox{\ie~~} \mathbf{G}(\mathbf{h}(\mathbf{X})) &=& D^2_{t\mathbf{X}^T}\mathbf{h}(\mathbf{X})
\cdot \mathbf{f}(\mathbf{X})+ D_{\mathbf{X}^T}\mathbf{h}(\mathbf{X})\cdot \mathbf{F}(\mathbf{X}).
\label{eq:dc1}
\end{eqnarray} 
Here $D^2_{t\mathbf{X}^T}\mathbf{h}(\mathbf{X})= (\mathbf{I}\otimes {{\dot{\phi}^T_t(\mathbf{X}})})\cdot 
D^2_{\mathbf{X}\mathbf{X}^T}\mathbf{h}(\mathbf{X})$ 
where $\otimes$ is direct product while $\mathbf{I}$ and $D^2_{\mathbf{X}\mathbf{X}^T}\mathbf{h}(\mathbf{X}) $ are identity and
 
Hessian matrices of
dimensions $n\times n$ and $nn\times n$ respectively \cite{ap,mag}.
If  $\mathbf{h}(\mathbf{X})$ is linear  then $D^2_{\mathbf{X}\mathbf{X}^T}\mathbf{h}(\mathbf{X})= \z$. Therefore,
 Eq. (\ref{eq:dc1}) can be recast as
\begin{eqnarray}
\mbox{} \mathbf{G}(\mathbf{h}(\mathbf{X})) &=& D_{\mathbf{X}^T}\mathbf{h}(\mathbf{X})\cdot \mathbf{F}(\mathbf{X}).
\label{eq:pp}
\end{eqnarray} 
 
\noindent
Since  $ D_{\mathbf{X}^T}\mathbf{h}(\mathbf{X})\ne \z$ for linear $\mathbf{h}(\mathbf{X})$, therefore the  perpetual
 points corresponding to $ \mathbf{F}(\mathbf{X})=0$ are mapped to $ \mathbf{G}(\mathbf{Y})=0$. This completes the proof.
 $\fbox{~}$\\

\noindent
{\bf Remark 1:}
If $\mathbf{h}(\mathbf{X})$ is nonlinear then the term $ D^2_{t\mathbf{X}^T}\mathbf{h}(\mathbf{X})
\cdot \mathbf{f}(\mathbf{X})$ may not be zero at perpetual points (where $\mathbf{f}(\mathbf{X})\ne 0$).
 If this term is not zero  then systems, $ \mathbf{F}(\mathbf{X})$ and $ \mathbf{G}(\mathbf{Y})$, are not 
topologically equivalent, and hence a  set of new perpetual points may be created for transformed 
system, $ \mathbf{G}(\mathbf{Y})$. 
Note that, for nonlinear $\mathbf{h}(\mathbf{X})$, a set of new perpetual points are also possible  in transformed system if
 $D^2_{t\mathbf{X}^T}\mathbf{h}(\mathbf{X})=\z$ and $ D_{\mathbf{X}^T}\mathbf{h}(\mathbf{X})=\z$ for $\mathbf{F}(\mathbf{X})\ne\z$.
However if $D^2_{t\mathbf{X}^T}\mathbf{h}(\mathbf{X})=\z, D_{\mathbf{X}^T}\mathbf{h}(\mathbf{X})=\z  $ and  $\mathbf{F}(\mathbf{X})= \z$
then the both $\mathbf{G}(\mathbf{Y})=\z$ and $\mathbf{g}(\mathbf{Y})=\z$, and hence the perpetual points
 corresponding to  $\mathbf{f}(\mathbf{X})=\z$ get mapped to the
new fixed points of $ \mathbf{g}(\mathbf{Y})$. 
 $\fbox{~}$\\

\noindent
{\bf Theorem 3:}  If the velocity vectors $\mathbf{f}(\mathbf{X})$ and $\mathbf{g}(\mathbf{Y})$ are topologically conjugate 
under diffeomorphic\footnote{ $\mathbf{h}(\mathbf{X})$ should be homeomorphic$^4$ as well as  
$\mathbf{h}$ and  $\mathbf{h}^{-1}$ have  continuous derivatives.} linear $\mathbf{h}$ then the eigenvalues
 $\lambda$ and $\mu$ for the respective fixed  and perpetual
 points, corresponding to the velocity and acceleration vectors, are preserved.\\

\noindent
{\sl Proof:} 
The derivative of Eq. (\ref{eq:fp})  with respect to  $\mathbf{X}$  gives
\begin{eqnarray}
\mbox{}  D_{\mathbf{h}^T}\mathbf{g}(\mathbf{h}(\mathbf{X}))\cdot D_{\mathbf{X}^T}\mathbf{h}(\mathbf{X})
&=& (\mathbf{I}\otimes {{\mathbf{f}^T}(\mathbf{X})})\cdot 
D^2_{\mathbf{X}\mathbf{X}^T}\mathbf{h}(\mathbf{X})+ D_{\mathbf{X}^T}\mathbf{h}(\mathbf{X})\cdot D_{\mathbf{X}^T}\mathbf{f}(\mathbf{X}).
\label{eq:lm}
\end{eqnarray}

\noindent
Here, $ D_{\mathbf{h}^T}\mathbf{g}(\mathbf{h}(\mathbf{X}))= D_{\mathbf{Y}^T}\mathbf{g}(\mathbf{Y})$. 
At fixed point $\mathbf{f}(\mathbf{X}_{FP})=\z$, hence
\begin{eqnarray}
 D_{\mathbf{Y}^T}\mathbf{g}(\mathbf{Y}) &=&
 D_{\mathbf{X}^T}\mathbf{h}(\mathbf{X})\cdot D_{\mathbf{X}^T}\mathbf{f}(\mathbf{X})\cdot (D_{\mathbf{X}^T}\mathbf{h}(\mathbf{X}))^{-1}.
\label{eq:dc}
\end{eqnarray}

\noindent
This shows that the set of eigenvalues $\lambda$ is same at the fixed points for the both flows.

\noindent
Similarly, the derivative of  Eq. (\ref{eq:pp}) with respect to  $\mathbf{X}$ gives
\begin{eqnarray}
\mbox{}  D_{\mathbf{Y}^T}\mathbf{G}(\mathbf{Y})\cdot D_{\mathbf{X}^T}\mathbf{h}(\mathbf{X})
&=&  (\mathbf{I}\otimes {{\mathbf{F}^T}(\mathbf{X})})\cdot 
D^2_{\mathbf{X}\mathbf{X}^T}\mathbf{h}(\mathbf{X})
+ D_{\mathbf{X}^T}\mathbf{h}(\mathbf{X})\cdot D_{\mathbf{X}^T}\mathbf{F}(\mathbf{X}).
\label{eq:dc0}
\end{eqnarray}

\noindent
At perpetual points $\mathbf{F}(\mathbf{X}_{PP})=\z$, hence
\begin{eqnarray}
 D_{\mathbf{Y}^T}\mathbf{G}(\mathbf{Y}) &=&
 D_{\mathbf{X}^T}\mathbf{h}(\mathbf{X})\cdot D_{\mathbf{X}^T}\mathbf{F}(\mathbf{X})
\cdot (D_{\mathbf{X}^T}\mathbf{h}(\mathbf{X}))^{-1}\mbox{~~~}
\label{eq:dcf}
\end{eqnarray}
\noindent
This shows that the set of eigenvalues $\mu$ is same at the perpetual  points for the both flows.
$\fbox{~}$\\

\noindent
{\bf Example 1:}
In order to demonstrate the mapping and eigenvalues  of/at perpetual points  
we  consider a simple, analytically tractable, dynamical system
\begin{eqnarray}
\dot{ x}&=&x^2-A^2,
\label{eq:1d}
\end{eqnarray}
where $A$ is a parameter.
It has  two fixed points, $x_{\tiny{ FP}}=\pm A$,
which are respectively stable and unstable with $\lambda=\pm2A$ \ie, initial conditions $x(0)<A$ lead to 
$x_{FP}=-A$ otherwise settle at  infinity.
 Its time derivative  gives acceleration,
\begin{eqnarray}
\ddot{ x}&=&2x(x^2-A^2),
\label{eq:1dac}
\end{eqnarray}
\noindent
which gives  one  perpetual point, $x_{PP}=0$.
The velocity at this perpetual points is maximum,  $A^2$, while
 the  eigenvalue is $\mu=-2A^2$ \cite{ap}.

Now let us transform this system using linear function, 
\begin{eqnarray}
y=h(x)&=&\al x+\beta
\label{eq:1dt}
\end{eqnarray}
where $\al$ and $\beta$ are parameters. The first and second order time derivatives  give, velocity and acceleration for
 the transformed system, as
\begin{eqnarray}
\dot{y}&=&\frac{(y-\beta)^2-\alpha^2 A^2}{\alpha} \\
 \mbox{and~~~~~}\ddot{y}&=&\frac{2}{\alpha^2}(y-\beta)[(y-\beta)^2-\alpha^2 A^2], 
\end{eqnarray}

\noindent
respectively, The fixed and the perpetual points of this transformed system are $y_{FP}=\beta\pm \alpha A$ and $y_{PP}=\beta$
respectively. Note that these are the transformed points corresponding to Eqs. (\ref{eq:1d}) and (\ref{eq:1dac}) under the  
transformation, Eq. (\ref{eq:1dt}).
The eigenvalues $\lambda=\pm2A$ and $\mu=-2A^2$ are the same as those of the Eqs. (\ref{eq:1d}) and (\ref{eq:1dac}) respectively. 
Therefore  under the linear transformation the  fixed and perpetual points are mapped and the eigenvalues are preserved 
[cf. Theorems 1, 2 and 3].
$\fbox{~}$ \\

\noindent
{\bf Example 2:}
Consider the same system as described by Eq. (\ref{eq:1d}). However, with nonlinear transformation, say,
\begin{eqnarray}
y=h(x)&=& x^2
\label{eq:trans2}
\end{eqnarray} 
 the velocity and acceleration for
 the transformed system become
\begin{eqnarray}
\dot{y}&=&2 \sqrt(y) (y- A^2) \\
 \mbox{and~~~~~}\ddot{y}&=& 2(3y-A^2)(y-A^2).
\end{eqnarray}
Note that the range of $y$ is $y\ge 0$. Its fixed points are $A^2$ and $0$ where the former is the 
mapping of fixed points of Eq. (\ref{eq:1d}) while
latter one is newly created. The eigenvalues, $\lambda$, for these fixed points are $\pm2A$ and $\pm \infty$ respectively
which depend on the sign of square root. The  perpetual point is $A^2/3$ which preimage is not the $x_{PP}=0$ of
 Eq. (\ref{eq:1d}), and its eigenvalues, 
$\mu$ is $-4A^2$. These suggest that due to nonlinear transformation new fixed or perpetual points
 points may be created [cf. Remarks 1 and 2].
$\fbox{~}$\\

\section{Summary}
In this work we established the theorems for the topological conjugacy of the perpetual points  under   linear transformation.
The eigenvalues also get preserved  under diffeomorphic linear transformation. However if the transformation
is nonlinear then this topological feature of the dynamical systems remains restricted.
Examples  corresponding to the  linear and nonlinear transformations are presented to confirm the conjugacy.

The concept of topological conjugacy is very important for comparing  dynamical systems which are
structurally different  (of various mathematical models) but have similar (equivalent) dynamics. This becomes 
even more important for nonlinear systems which are not solvable analytically. For such systems
global dynamical behavior  are normally guessed/analyzed from its linearized version,  \eg,
understanding the global dynamical behavior of the full system is determined by assembling the nearby (local) motions 
around all the  fixed points present in the system. Note that, as the motion at fixed points is stationary while 
extremum (or inflection type) velocity at perpetual points is nonzero, 
 the understanding of the dynamics of oscillating systems near perpetual points
are necessary to understand the full system. Therefore, the results of topological conjugacy presented in this paper are 
of immediate use  for understanding the global dynamics of the nonlinear systems.

It may be mentioned that many problems/theorems, similar to those studied in the context of  fixed points, 
now need extension for the studies of perpetual points \eg~ Hartman-Grobman theorem, $C^k$-conjugacy \etc

\section*{Acknowledgments}
Author  thanks, R. S. Kaushal, M. D. Shrimali, B. Biswal, N. Kuznetsov, and R. Ramaswamy  for critical comments on this
 manuscript.  The financial supports from the DST, Govt. of India and Delhi University  Research \& Development Grant
are gratefully  acknowledged.\\

 \bibliographystyle{IJNSBST} 

\clearpage

\end{document}